\documentclass{article} 
\usepackage[utf8]{inputenc}         
\usepackage{amssymb,amsfonts,amsmath,eqnarray} 
\usepackage{graphicx}                 
\usepackage{enumerate}                %
\usepackage{hyperref}
\usepackage{tikz}
\usepackage{xcolor}

\usepackage[nomessages]{fp}

{}

\title{Knapsack problem for nilpotent groups}

\author{Alexei Mishchenko, Alexander Treier}

\newtheorem{theorem}{Theorem}
\newtheorem{lemma}{Lemma}
\newtheorem{proposition}{Proposition}
\newtheorem{corollary}{Corollary}
\newtheorem{definition}{Definition}

\newcommand{\el}{l}
\newcommand{\dconst}{\blacksquare}
\newcommand{\rank}{{\mathrm{rank}}}
\newcommand{\KP}{{\mathrm{KP}}}
\newcommand{\Z}{{\mathbb{Z}}}
\newcommand{\N}{{\mathbb{N}}}
\newcommand{\K}{{\mathcal{K}}}
\newcommand{\ut}[1]{UT_{#1}(\mathbb{Z})}
\newcommand{\oneTo}[1]{1, \ldots, {#1}}
\newcommand{\listSym}[2]{{{#1}}_{1}, \ldots, {{#1}}_{{#2}}}
\newcommand{\expForm}[3]{{#1}^{{#2}_1}_{1} \ldots {#1}^{{#2}_{#3}}_{#3}}
\newcommand{\expFormKP}[1]{\expForm{g}{\varepsilon}{#1}}
\newcommand{\refb}[1]{(\ref{#1})}
\newcommand{\numCommutators}{322}
\newcommand{\numGenerators}{26}

\def\proof{{\noindent{\bf Proof.}} }
\allowdisplaybreaks[1]

\begin{document}
\maketitle

\tableofcontents

\section{Introduction}

In the paper~\cite{MNU:KnapsackProblemInGroups} A.~Myasnikov, A.~Nikolaev, and A.~Ushakov stated a group version of the well known Knapsack problem. The motivation for our research and initial results in this direction may be found in that paper, and further results in \cite{MNU:KnapsackProblemInGroups2, Lor:KPGg, Lor:KPNilp}.

We give a definition of Knapsack problem for groups following \cite{MNU:KnapsackProblemInGroups}. Let $G$ be an arbitrary group with a presentation $G = \langle X | R \rangle$ and solvable word problem. Let $g_1, \ldots, g_k, g$ be finite words in the alphabet $X \cup X^{-1}$. Then the Knapsack Problem for the group $G$ is stated in the following way.

\noindent
\textbf{Knapsack Problem. $\KP$.}
\textit {Given input words $\listSym{g}{k}, g$, decide whether there exist integers $\listSym{\varepsilon}{k}$ such that the equality
\begin{equation}\label{eq:KnapsackProblem}
\expFormKP{k} = g
\end{equation}
holds in the group $G$.}

There are several notable questions related to $\KP$. One such question is that of decidability of $\KP$ for a specific class of groups $\K$. In the case when $\KP$ is decidable for a class $\K$, another natural question is how computationally hard $\KP$ for class $\K$ is. In this regard, it is known that $\KP$ is decidable in polynomial time for abelian and hyperbolic groups. In this work we investigate decidability of $\KP$ for nilpotent groups.

The main results of the present paper are as follows. In \textit{Theorem \ref{theorem:main}} we prove that Knapsack problem ($\KP$) is undecidable for any group of nilpotency class two if the number of generators (without torsion) of the derived subgroup is at least $322$. This theorem together with the fact that if $\KP$ is undecidable for a subgroup then it undecidable for the whole group allows us extend our result to certain classes of polycyclic groups, linear groups and nilpotent groups of higher nilpotency class ($\geq 3$).

We draw the reader's attention to a result of Daniel König, Markus Lohrey, and Georg Zetzsche~\cite{Lor:KPNilp} that $\KP$ is undecidable for a direct product of sufficiently many copies of the discrete Heisenberg group $H_3(\Z)$. This implies that $\KP$ is generally undecidable for nilpotent groups. We would like to point out that our approach is different from that of Daniel König, Markus Lohrey, and Georg Zetzsche. Moreover, our Theorem 1 provides an explicit bound, $322$, for the number of copies of $H_3(\Z)$ in a direct product that suffices for undecidable $\KP$. The paper~\cite{Lor:KPNilp} also contains interesting results on Subset Sum Problem and Knapsack problem for nilpotent, polycyclic, and co-context-free groups.

The authors are grateful to A.~Miasnikov and A.~Nikolaev for their advice and discussions.

\section{Preliminaries}

\subsection{Nilpotent groups}

Recall the definition and basic properties of nilpotent groups. A group $G$ is called a nilpotent group of class $c$ if it has a lower central series of length $c$:
$$G = G_1 \trianglerighteq G_2 \trianglerighteq \ldots \trianglerighteq G_c \trianglerighteq G_{c+1} = \{1\},$$
where $G_{k+1} = [G_k, G]$, $k = \oneTo{c}$ and $G_1 = G$. 

Let $X = \{\listSym{x}{n}\}$ be a set of letters, and $G = \langle X \rangle$ be a free nilpotent group of class $2$. By definition, the following identity holds for group $G$: 

\begin{equation}\label{eq:Identity}
\forall  x,y,z \in G \ [x,[y,z]] = 1 
\end{equation}

Using identity \refb{eq:Identity}, the collection process in group $G$ is organized via the transformation  
\begin{equation}\label{eq:Collecting}
yx = xy [x,y]^{-1},
\end{equation}
where $x,y$ are any elements of $G$. Using the equality~\refb{eq:Collecting} we can reduce any word $g$ in the alphabet $X \cup X^{-1}$ to the normal form for elements of the group~$G$: 
\begin{equation}\label{eq:NormalForm}
 g = \expForm{x}{\alpha}{n} \prod_{i<j} [x_i,x_j]^{\beta_{ij}},
\end{equation}
where $\alpha_i,\beta_{ij} \in \Z, i,j = \oneTo{n}, i < j$ and $[x_i,x_j] = x^{-1}_ix^{-1}_jx_ix_j$.

Using \refb{eq:Identity}, it is not hard to show that for any two elements $a,b$ of the group~$G$ and $\alpha, \beta \in \Z$ we have the following equality:

\begin{equation}\label{eq:RingMultiplication}
[a^\alpha,b^{\beta}] = [a,b]^{\alpha\beta}.
\end{equation}

\subsection{Knapsack problem}\label{section:KP}

We stated the Knapsack problem ($\KP$) for groups in Introduction. Recall that the $\KP$ is called decidable for the class of groups $\K$ if for any group $G \in \K$ there exists an algorithm that, given any input $\listSym{g}{k}, g$, answers the question whether or not the exponential group equation~\refb{eq:KnapsackProblem} has a solution in the group $G$. We can restrict the notion of decidability of $\KP$ and explore $\KP$ for single group or for some type of inputs of $\KP$. In our work we concentrate on decidability of $\KP$ for the class of nilpotent groups.

Let $G$ be a free nilpotent group of class $2$ and let $\listSym{g}{k}, g$ be presented in the form~\refb{eq:NormalForm}.
Using~\refb{eq:Collecting} and~\refb{eq:RingMultiplication} we can reduce the expression $\expFormKP{k}$ to the form \refb{eq:NormalForm}. Thus, the following proposition holds:

\begin{proposition}\label{prop:KPtoDiophantine}
Let $G$ be a free two-step nilpotent $n$-generated group. Then $\KP$ stated above for the group $G$ is equivalent to a system of Diophantine equations with unknowns $\listSym{\varepsilon}{k}$ of degree $2$. Moreover, the number of linear equations in the system is not greater than $n$ and the number of quadratic equations is not greater than $\frac{n(n-1)}{2}$.
\end{proposition}

\subsection{Diophantine equations and Hilbert's Tenth problem}\label{section:DiophantineEquations}

\textit{Proposition \ref{prop:KPtoDiophantine}} shows that $\KP$ for nilpotent groups is closely related to Diophantine equations. This section is devoted to Diophantine equations.  

A polynomial equation $D(\listSym{x}{n}) = 0$ with integer coefficients is called Diophantine.

In 1900 at the Second International Congress of Mathematicians D.~Hilbert presented his famous list of problems. The 10th problem is concerned with Diophantine equations. The problem statement is as follows: is there an algorithm that for any Diophantine equation answers the question whether or not this equation has a solution in integers? In 60-70th of previous century M.~Davis, J.~Robinson, H.~Putnam, and Yu.~Matyasevich proved that there is no algorithm to decide whether an arbitrary Diophantine equation has solution in integers or not. For more details on Hilbert's Tenth Problem we refer the reader to the book of Yu.~Matiyasevich~\cite{Mat:Hilbert10}, which, in addition to the solution of the problem, provides a historical survey and describes a number of applications of negative solution of Hilbert's Tenth Problem.

In some cases of Diophantine equations there exists an algorithm to decide whether the equation has a solution. In \cite{Siegel:QuadraticDecidable} C.~Siegel gives an algorithm for a single Diophantine equation of degree $\leq 2$. So, if we have $2$-generated free two-step nilpotent group $G$ (which is known as Heisenberg group) by \textit{Proposition \ref{prop:KPtoDiophantine}} the $\KP$ for any input is equivalent to a system of two linear equations and one quadratic equation. Such a system may be reduced to a single quadratic equation (for example, this is shown in~\cite{DLS:EquationsInNilpotentGroups}), and therefore, the following proposition holds:

\begin{proposition}\label{prop:Heisenberg group}
The Knapsack problem for Heisenberg group is decidable on any input.
\end{proposition}

Now we return to the question of undecidability of Diophantine equations. From papers of Julia Robinson, Martin Davis, Hilary Putnam \cite{DPR:ExpDiophantine} and Yu.~Matiyasevich \cite{ Mat:EnumerableAreDiophantine} every recursive enumerable set $W$ can be presented in Diophantine form:
\begin{equation}
x \in W \iff \exists \listSym{x}{n}\ P(x, \listSym{x}{n}) = 0,
\end{equation}
where the variables $\listSym{x}{n}$ are positive integers and $P(x, \listSym{x}{n})$ is a Diophantine polynomial. Since there exist recursively enumerable but non recursive sets then there is no algorithm to decide for arbitrary Diophantine equation whether it has a solution. Moreover, if $W_1, W_2, \ldots$ is a list all recursively enumerable sets, then there is a polynomial $U$ such that for any $k \in \N$
\begin{equation}\label{eq:UniversalEquation}
x \in W_k \iff \exists \listSym{x}{n}\ U(x, k, \listSym{x}{n}) = 0.
\end{equation}

The polynomial $U(x,k,\listSym{x}{n})$ has fixed degree and fixed number of variables. Such polynomial $U$ is called a universal polynomial. J.P.~Jones in~\cite{Jon:UndecidableDEquations, Jon:UniversalDEquations} constructed a universal system of equations that can be reduced to a universal polynomial of degree $4$ with $58$ unknowns. To reduce the Jones system to a single equation we need to prepare this system (because some equations have degree greater than $2$) by transformations and substitutions which are described by Jones. After that we introduce several new variables which are tied by linear relations to lower the number of generators of two step nilpotent group $G$ for building an input for $\KP$ (see the next sections for details). We are not aware of any published work that provides an explicit version of the universal system of equations of degree $\leq 2$, so we give this system in the present paper. In the next sections we use this system for constructing a universal $\KP$ input and calculating rank of nilpotent groups with undecidable $\KP$.

Any letter symbols in system below are variables except $x, \dconst_z, \dconst_y, \dconst_u$ which are positive integer parameters of $U$. The constants $\dconst_z, \dconst_y, \dconst_u$ encode a r.e. set which determines the universal system. So if we put $\dconst_z, \dconst_y, \dconst_u$ that encode a non-recursive set $W$, then there is no algorithm for any $x \in W$ to answer the question whether the equation have a solution. After applying transformations to Jones system we obtain the following universal system:
\begin{gather}
\label{system:first}
\Gamma_1 = \Gamma_{26}^2,
\\
\Gamma_2 = MU,
\\
\Gamma_3 = B (2\Gamma_{23} - B) - 1,
\\
\Gamma_4 = \Gamma_{23} C_1,
\\
\Gamma_5 = c^2 ,
\\
\Gamma_6 = \Gamma^2_5, 
\\
\Gamma_8 = \Gamma_{24}^2,
\\
\Gamma_9 = \lambda B,
\\
\Gamma_{10} = GH,
\\
\Gamma_{11} = F^2,
\\
%
\Gamma_{12} = \Gamma_{23} E,
\\
\Gamma_{13} = \Gamma_{25}^2,
\\
\Gamma_{14} = \Gamma_{23}\Gamma_{25},
\\
\Gamma_{15} = N^2,
\\
\Gamma_{16} = YK,
\\
\Gamma_{18} = PK,
\\
\Gamma_{20} = \Gamma_8 \Gamma_{24},
\\
\Gamma_{21} = \Gamma^{2}_8,
\\
\Gamma_{22} = \Gamma_6\Gamma_{20},
\\
%
B = 2\Gamma^2_1(2\dconst_z)^{{5}^{59} + 1},
\\
D_1 = 1 + \Gamma_{27} + C_1(\Gamma_{23} - B) + \alpha\Gamma_3,
\\
(\Gamma_4 - C_1)(\Gamma_4 + C_1) + 1 = D^2_1,
\\
C_1 = 5^{59} + \Delta (\Gamma_{23} - 1),
\\
c = 1 + (\Gamma_{26} - \varepsilon)B + g,
\\
%
e + 2\dconst_z \Gamma_{26}\el + 2\dconst_z B\Gamma_6 + \Gamma_7 = 2\dconst_z (1 + \Gamma_{27}),
\\
\el = \dconst_u + t(B-2\dconst_z),
\\
e = \dconst_y + m(B - 2\dconst_z),
\\
%
%
%
%
%
S = g - 4 \dconst^2_z \Gamma_{22} + \el \Gamma_{24} + e(\Gamma_8 + 4 \dconst_z \Gamma_{22}) + 2 \dconst_z \Gamma_9 (-2\dconst_z \Gamma_{22} + \Gamma_{20} + \Gamma_{21}), 
\\
T = \Gamma_{24} - 1 - (\Gamma_{26} - 1) \el + (\Gamma_9 - 2\lambda\dconst_z)(\Gamma_{24} + \Gamma_8)  + 2\dconst_z(B-2)  \Gamma_{21},
\\
N = 16 \dconst_z \Gamma_{20}\Gamma_{8},
\\
R = S(\Gamma_{15} - N) + (T+1)(\Gamma_{15} - 1),
\\
P = 2 M \Gamma_2,
\\
\label{system:example}
(K - \Gamma_{18})(K + \Gamma_{18}) + \Gamma^2_{19} = 1,
\\
(2\Gamma_{25} - 2 \Gamma_{16} - K)(2\Gamma_{25} - 2 \Gamma_{16} + K) + \Gamma_{17}  = 0,
\\
K = R + 1 + h(P-1),
\\
M = RY,
\\
%
%
U = \Gamma_{15}w,
\\
Y = \Gamma_{15}s,
\\
%
%
D = -2\Gamma_{25} - 5\gamma + \Gamma_{26}w + \Gamma_{23}(\Gamma_{25} + 4\gamma),
\\
I = D +oF,
\\
(D - \Gamma_{14})(D + \Gamma_{14}) + \Gamma_{13} = 1,
\\
E = i\Gamma_{13} + 1,
\\
(\Gamma_{12} - E)(\Gamma_{12} + E) - \Gamma_{11} + 1 = 0,
\\
G = \Gamma_{23} + \Gamma_{11}(\Gamma_{11} - \Gamma_{23}),
\\
H = 2R + 1 + j\Gamma_{25},
\\
I^2 + H(H - \Gamma_{10}) = 1, \\
\Gamma_{23} = \Gamma_2 + M, \\
\Gamma_{24} = 1 + \Gamma_9 - \lambda, \\  
\Gamma_{25} = 2R + 1 + C_1 + \varphi, \\
\Gamma_{26} = \varepsilon + x, \\
\label{system:last}
\Gamma_{27} = \lambda (B-1).
\end{gather}

\section{Equivalence between system of Diophantine equations and Knapsack Problem for nilpotent groups}

In this section we show that any finite system of Diophantine equations is equivalent to $\KP$ for some two step nilpotent group $G$ on some input. This means that for any finite system $S$ of Diophantine equations there exists a group $G = \langle \listSym{x}{n} \rangle$ and input $\listSym{g}{k}, g$ which are words of alphabet $X \cup X^{-1}$ such that $\KP$ for group $G$ has solution if and only if the system $S$ has solution. 

Let $S = \{\listSym{s}{r}\}$ be a finite system of Diophantine equations with variables $\listSym{x}{n}$, where $s_{i} := (f_i(\listSym{x}{n}) = c_i) $ is a Diophantine equation. Since any finite system of Diophantine equations is equivalent to finite system of equations of degree less or equal than $2$, we may assume that every equation in $S$  written in the form  
\begin{equation}\label{eqn:quadratic}
s_i := \left(\sum^{n}_{i=1} \alpha_ix_i + \sum^{n}_{i,j = 1} \beta_{ij}x_ix_j = \gamma\right),
\end{equation}
where $\alpha_i,\beta_{ij}, \gamma \in \Z $. 

We start by showing how to construct an input for $\KP$ equivalent to a single quadratic Diophantine equation~\refb{eqn:quadratic}. Let $a,b$ be generators of the group $G$ and $[a,b]$ a nontrivial basic commutator in $G$. Below we pick elements $\listSym{g}{r}\in G$ such that the $\expFormKP{r}$ is equal to $[a,b]^{\sum^{n}_{i=1} \alpha_ix_i + \sum^{n}_{i,j = 1} \beta_{ij}x_ix_j}$, then we put $g = [a,b]^\gamma$. $\KP$ on the obtained input will be equivalent to \refb{eqn:quadratic}. 

Consider the linear part of~\refb{eqn:quadratic}. For every summand $\alpha_ix_i, \ i = \oneTo{n}$ we put $g_i = [a,b]^{\alpha_i}$ and get $\expFormKP{n} = [a,b]^{\sum^{n}_{i=1} \alpha_i\varepsilon_i}$. Thus, we assume that $x_i = \varepsilon_i$. 

Turn to the quadratic part of~\refb{eqn:quadratic}. For every summand $\beta_{ij}x_ix_j$ we assign four new elements of input (we assume that in previous steps we constructed $r$ elements of input):
\begin{eqnarray*}
g_{r+1} & = & a^{-\beta_{ij}} \cdot c_1, \\
g_{r+2} & = & b^{-1} \cdot c_2, \\
g_{r+3} & = & a^{\beta_{ij}}\cdot  c^{-1}_1, \\
g_{r+4} & = & b \cdot c^{-1}_2,
\end{eqnarray*}
where $c_1,c_2 \in [G,G]$ are non-trivial commutators that have not appeared previously in construction of the input. Then $$K = g^{\varepsilon_{r+1}}_{r+1} g^{\varepsilon_{r+2}}_{r+2} g^{\varepsilon_{r+3}}_{r+3} g^{\varepsilon_{r+4}}_{r+4} = a^{-\beta_{ij}\varepsilon_{r+1}}b^{-\varepsilon_{r+2}}a^{\beta_{ij}\varepsilon_{r+3}}b^{\varepsilon_{r+4}} c^{\varepsilon_{r+1} - \varepsilon_{r+3}}_1 c^{\varepsilon_{r+2} - \varepsilon_{r+4}}_2.$$ 

Setting that the exponents of commutators $c_1$ and $c_2$ are equal to zero in element $g$ is equivalent to the condition $\varepsilon_{r+1} = \varepsilon_{r+3}$ and $\varepsilon_{r+2} = \varepsilon_{r+4}$. As a result we have $K = [a,b]^{\beta_{ij}\varepsilon_{r+1}\varepsilon_{r+2}}$. Now we need to tie the values of $\varepsilon_{r+1}$ to $\varepsilon_i$ and $\varepsilon_{r+2}$ to $\varepsilon_j$. To do that we apply the same trick as in the previous case. Let $c_3$ be a non-trivial commutator that we have never used before. We put $g'_i = g_i c_3$ and $g'_{r+1} = g_{r+1} c^{-1}_3$, then we replace $g_{i}$  by $g'_{i}$ and $g_{r+1}$ by $g'_{r+1}$ in the input. The imposed restrictions give us $\varepsilon_{r+1} = \varepsilon_i = x_i$. Then we repeat the same with $\varepsilon_{r+2}$ and $\varepsilon_j$. Proceeding in the same way with all other quadratic summands we finally get the following exponential expression:
$$
\expFormKP{k} = [a,b]^{\sum^{n}_{i=1} \alpha_ix_i + \sum^{n}_{i,j = 1} \beta_{ij}x_ix_j},
$$
where $x_i = \varepsilon_i, \ i = \oneTo{n}$. Then we set $g = [a,b]^\gamma$ and obtain the exponential equation $\expFormKP{k} = g$ in group $G$ equivalent to Diophantine equation \refb{eqn:quadratic}. 

It is easy to see that if we have an arbitrary finite system $S$ of $l$ quadratic Diophantine equations we can build an input for $\KP$ that realizes all equations in the system $S$ as powers of $l$ basic commutators ($[a,b], [a,c], [c,d]$, e.t.c., where $a,b,c,d, \ldots$ are generators of $G$) as described above. Thus, for any finite system $S$ and any nilpotent group $G$ with sufficiently many basic commutators (recall that, besides $l$ basic commutators for equations of $S$, we need more commutators to realize bindings between variables of $\KP$) we can construct an input on which $\KP$ for the group $G$ is equivalent to system $S$. 

From the above we have the following
\begin{proposition}\label{prop:equivalence}
For any finite system of Diophantine equations exists a finitely generated free group $G$ of nilpotency class $2$ and an input $\listSym{g}{k}, g \in G$ such that $\KP$ on this input has solution in $G$ if and only if the system $S$ has solution in $\N \cup \{0\}$.
\end{proposition}

Now we briefly describe another, more general, approach to establishing equivalence between $\KP$ for nilpotent groups and decidability of Diophantine equations. This reduction may be more convenient than the one described above in case of an arbitrary Diophantine equation (or any finite system of equations) of degree greater than~$2$.

We begin by defining the notion of a Diophantine term by induction as follows.
\begin{definition}
\begin{enumerate}
\item Every constant is a term.
\item Every variable is a term.
\item For every two terms $t_1$ and $t_2$ the $t_1 + t_2$ and $t_1 t_2$ are terms.
\end{enumerate}
A term is called {\em simple} if it is a constant or a variable.
\end{definition}

We can present any Diophantine equation as equality of two terms $t_1 = t_2$. There are many ways to express a given polynomial as a combination of sums and products of Diophantine terms. For example, we may present the polynomial $f(x) = x^2 - 1$ as sum of two terms: $x^2$ and $-1$ and then $x^2$ is a product of $x$ and $x$, or we may look at $f(x)$ as the product of $x-1$ and $x+1$. We can represent computation scheme of a term as a binary tree where leafs are simple terms and internal vertices are symbols of multiplication~``$\cdot$'' or addition~``$+$''. 

Let $t_1 = t_1(\listSym{\varepsilon}{n})$ and $t_2 = t_2(\listSym{\varepsilon}{n})$ be Diophantine terms such that $\expFormKP{r} = h \cdot [a,b]^{t_1}[c,d]^{t_2}$ and the powers of $[a,b]$ and $[c,d]$ in $g$ and $h$ are equal to zero.
Thus, to describe how to construct an input for $\KP$ equivalent to a given Diophantine polynomial we need to show, for two terms $t_1, t_2$, how to extend the input to realize the following: terms $t_1 + t_2$, $t_1 \cdot t_2$ and equations $t_1 = t_2$, $t_1 = \gamma$, where $\gamma \in \Z$. 

\begin{itemize}
\item[$(t_1 = \gamma)$]: to satisfy this condition we introduce one new input element:
\begin{eqnarray*}
g_{r+1} & = & [a,b] c_1,
\end{eqnarray*}
where $c_1$ is a basic commutator in $G$ which has not been used before, and set $g' = g c^{\gamma}_1$.

\item[$(t_1 = t_2)$]: in this case we introduce two new input elements:
\begin{eqnarray*}
g_{r+1} & = & [a,b]^{-1} c_1, \\
g_{r+2} & = & [c,d]^{-1} c^{-1}_1, 
\end{eqnarray*}
then $\expFormKP{r+2} = [a,b]^{t_1 - \varepsilon_{r+1}}[c,d]^{t_2 - \varepsilon_{r+2}} c^{\varepsilon_{r+1} - \varepsilon_{r+2}}_1$, which gives us $t_1 = \varepsilon_{r+1} = \varepsilon_{r+2} = t_2$ (provided that the powers of $[a,b], [c,d]$, and $c_1$ in $g$ are~$0$ ). 

\item[$(t_1 + t_2)$]: 
\begin{eqnarray*}
g_{r+1} & = & [a,b]^{-1} c_1, \\
g_{r+2} & = & [c,d]^{-1} c_1, 
\end{eqnarray*}
then $\expFormKP{r+2} = [a,b]^{t_1 - \varepsilon_{r+1}}[c,d]^{t_2 - \varepsilon_{r+2}} c^{\varepsilon_{r+1} + \varepsilon_{r+2}}_1$, which gives us $c^{t_1 + t_2}_1$ provided that the powers of $[a,b]$ and $[c,d]$ in the element $g$ are $0$.

\item[$(t_1 \cdot t_2)$]:
\begin{eqnarray*}
g_{r+1} & = & [a,b] x^{-1} \cdot c_1, \\
g_{r+2} & = & [c,d] y^{-1} \cdot c_2, \\
g_{r+3} & = & x \cdot  c^{-1}_1, \\
g_{r+4} & = & y \cdot c^{-1}_2,
\end{eqnarray*}
then 
\begin{eqnarray*}
\expFormKP{r+2} = 
& [a,b]^{t_1 - \varepsilon_{r+1}}[c,d]^{t_2 - \varepsilon_{r+2}}    \\
&  c^{\varepsilon_{r+1} - \varepsilon_{r+3}}_1 c^{\varepsilon_{r+2} - \varepsilon_{r+4}}_2 \\
&  x^{-\varepsilon_{r+1}}y^{-\varepsilon_{r+2}}x^{\varepsilon_{r+3}}y^{\varepsilon_{r+4}}
\end{eqnarray*}
which gives us $[x,y]^{t_1 \cdot t_2}$ provided that the powers of $[a,b], [c,d]$, $c_1,c_2$ in $g$ are $0$. 
\end{itemize}

\section{Nilpotent groups with undecidable $\KP$}

In previous section we described two reductions of any Diophantine equation or a system of Diophantine equations to $\KP$ in a nilpotent group $G$ with sufficient number of generators. Now we want to give a lower bound for the number of basic commutators in $G'$ of a torsion free two step nilpotent group $G$ with undecidable $\KP$. We do not aim to get the lowest possible bound for the number of commutators in a group $G$, but we note some simple transformations of the original Jones system of equations to reduce the number of generators. We omit a full description of the input for $\KP$ (because it contains $334$ elements of input) which is equivalent to the system of equations \refb{system:first} -- \refb{system:last}. However, we give an example that clarifies the process of input construction.

Consider an equation \refb{system:example}:
$$
(K - \Gamma_{18})(K + \Gamma_{18}) + \Gamma^2_{19} = 1.
$$
Let a,b be generators of $G$ such that the commutator $[a,b]$, along with commutators $\listSym{c}{7}$, have never been used before. Then we put
\begin{eqnarray*}
 g_{1} & = & a^{-1} c_1 c_3, (\text{for }  K)\\
 g_{2} & = & a^{-1} c_1 c_4, (\text{for } -\Gamma_{18}) \\
 g_{3} & = & b^{-1} c_2 c^{-1}_3, (\text{for } K)\\
 g_{4} & = & b^{-1} c_2 c_4, (\text{for } \Gamma_{18})\\
 g_{5} & = & a c^{-1}_1,\\
 g_{6} & = & b c^{-1}_2.
\end{eqnarray*}
Thus, the elements $\listSym{g}{6}$ are used to construct the term that corresponds to $(K - \Gamma_{18})(K + \Gamma_{18})$.
The next input elements $g_7, g_8, g_9, g_{10}$ serve in a similar capacity for $\Gamma^{2}_{19}$,
\begin{eqnarray*}
 g_{7} & = & a^{-1} c_5 c_7,\\
 g_{8} & = & b^{-1} c_6 c^{-1}_7,\\
 g_{9} & = & a c^{-1}_5,\\
 g_{10} & = & b c^{-1}_6.
\end{eqnarray*} 
Finally, the right hand side of $\KP$ expression is given by
$$
g  =  [a,b]^{1},
$$
and all commutators $\listSym{c}{7}$ have zero power in the element $g$.

\begin{eqnarray*}
\expFormKP{6} & = & a^{-\varepsilon_{1} - \varepsilon_{2}}b^{-\varepsilon_{3} - \varepsilon_{4}}a^{\varepsilon_{5}}b^{\varepsilon_{6}} c^{\varepsilon_1 + \varepsilon_2 - \varepsilon_5}_1 c^{\varepsilon_3+ \varepsilon_4 - \varepsilon_6}_2 c^{\varepsilon_1 - \varepsilon_3}_3 c^{\varepsilon_2 + \varepsilon_4}_4 = \\
& = & [a,b]^{(\varepsilon_1 + \varepsilon_2)(\varepsilon_3+\varepsilon_4)}c^{\varepsilon_1 - \varepsilon_3}_3 c^{\varepsilon_2 + \varepsilon_4}_4 = \\
& = & [a,b]^{(\varepsilon_1 + \varepsilon_2)(\varepsilon_1 - \varepsilon_2)} = (\text{put } \varepsilon_1 = K, \ \varepsilon_2 = \Gamma_{18} )\\
& = & [a,b]^{(K - \Gamma_{18})(K + \Gamma_{18})}.
\end{eqnarray*}

\begin{eqnarray*}
g^{\varepsilon_7}_7 g^{\varepsilon_8}_8 g^{\varepsilon_9}_9 g^{\varepsilon_{10}}_{10} & = & a^{-\varepsilon_{7}}b^{-\varepsilon_{8}}a^{\varepsilon_{9}}b^{\varepsilon_{10}} c^{\varepsilon_7 - \varepsilon_9}_5 c^{\varepsilon_8 - \varepsilon_{10}}_6 c^{\varepsilon_7-\varepsilon_8}_7 = \\
& = & [a,b]^{\varepsilon_7\varepsilon_8} c^{\varepsilon_7 - \varepsilon_8}_7 = \\
& = & [a,b]^{\varepsilon^2_7} = (\text{put } \varepsilon_7 = \Gamma_{19})\\
& = & [a,b]^{\Gamma^2_{19}}.
\end{eqnarray*}
The two latter expressions are equivalent to the following system:

$$\left\{\begin{array}{l}\varepsilon_1 + \varepsilon_2 = \varepsilon_5; \\
\varepsilon_3 + \varepsilon_4 = \varepsilon_6; \\
\varepsilon_1 = \varepsilon_3; \\
\varepsilon_2 = -\varepsilon_4; \\
\varepsilon_7 = \varepsilon_8 = \varepsilon_9 = \varepsilon_{10}; \\
(\varepsilon_1 + \varepsilon_2) (\varepsilon_1 - \varepsilon_2) + \varepsilon_7^2 = 1;
\end{array}\right.$$
Combining everything together we get
$$
\expFormKP{10} = [a,b]^{(K - \Gamma_{18})(K + \Gamma_{18}) + \Gamma^2_{19}} = [a,b],
$$
which gives us the desired equation~\refb{system:example}.

Finally, we need $167$ basic commutators in the group $G$ to interpret all equations~\refb{system:first}--\refb{system:last}. If any variable occurs $n+1$ times in our system, then we need another $n$ commutators to tie these variables. Additionally, we need $155$ commutators to tie the same variables in the equations. Hence the total number of commutators to realize the system~\refb{system:first}--\refb{system:last} is $167 + 155 = 322$. The input for $\KP$ is given by elements $g_1, \ldots, g_{334},g$, which depend on four integer parameters $x, \dconst_z, \dconst_y, \dconst_u$.




Based on previous computations we have the following

\begin{lemma}\label{lemma:main} Let $G$ be a torsion free group of nilpotency class $2$ with $\rank([G, G]) > \numCommutators$, then for every recursively enumerable set $W$ exists an input $I_W = \{g_1, \ldots, g_{334}, g\}$ such that
$$
x \in W \text{ iff } \KP \text{ has a solution in the group } G \text{ for the input } I_W.
$$
\end{lemma}
\proof For every recursively enumerable set $W$ there exist parameters $\dconst_z, \dconst_y, \dconst_u$ such that an integer $x$ lies in $W$ if and only if the system $S_W(x,\dconst_z, \dconst_y, \dconst_u)$ has a solution. Since $\rank([G, G]) > \numCommutators$ we can construct an input $I_W = \{g_1, \ldots, g_{334}, g\}$ for $\KP$ such that the corresponding instance of $\KP$ for $G$ has a solution if and only if the system $S_W$ has a solution. $\square$

\begin{theorem}\label{theorem:main}
Let $G$ be a torsion free group of nilpotency class $2$ and $\rank([G, G]) > \numCommutators$, then group $G$ has undecidable $\KP$ problem.
\end{theorem}

\proof There is set a $W$ that is recursively enumerable but is not enumerable. The statement follows by applying Lemma~\ref{lemma:main} to this set $W$. $\square$

\section{Corollaries}

In this section we give corollaries of \textit{Theorem~\ref{theorem:main}}. 

\begin{corollary}\label{corollary:freenilp}
Let $G$ be a free group of nilpotency class $2$ with $n$ generators. If $n$ is at least $\numGenerators$ then the group $G$ has undecidable $\KP$.
\end{corollary}

\proof Note that $G$ has $\frac{n(n-1)}{2}$ basic commutators. Since it is enough to have $\numCommutators$ basic commutators, we see that $\numGenerators$ generators suffice. $\square$

\begin{corollary}
Let $G$ be a group of nilpotency class $2$, $H$ be its torsion subgroup, $G_1 = G / H$ be the corresponding quotient group. If $\rank([G_1, G_1]) > \numCommutators$, then the group $G$ has undecidable $\KP$.
\end{corollary}

\begin{corollary}
If $n \geq 53$ then $\KP$ is undecidable for groups $\ut{n}$, $GL_n(\Z)$, $SL_n(\Z)$. 
\end{corollary}

\proof Denote by $F_{k}$ the free 2-step nilpotent group of rank $k$ with generators $X = \{\listSym{x}{k}\}$. By Corollary~\ref{corollary:freenilp} the $KP$ is undecidable for the group $F_{\numGenerators}$. By the theorem of Jennings every finitely generated torsion-free nilpotent group can be embedded into $\ut{n}$. Willem A. De Graaf and Werner Nickel~\cite{GN:ConstructingRepresentation} give the algorithm that constructs this embedding. Hence the $\KP$ is undecidable for $\ut{n}$ and we only need to get an estimate of $n$. The algorithm described by De Graaf and Nickel embeds the group $F_k$ in $\ut{n}$, where $n = k + C^2_k$. We construct an embedding $\rho$ which embeds $F_n$ into $\ut{2n+1}$. 

For every generator $x_i$ of the group $F_n$ we define an $(n+1)\times (n+1)$ matrix $M_i$,
$$
M_i = 
\bordermatrix{
 ~&       &   &        & i+1 &        & \cr
  &0      &   &        &     &        & \cr
  &\vdots &   &        &     &        & \cr
i &1      &   &        &     &        & \cr
  &\vdots &   &        &     &        & \cr 
  &0      &   &        &     &        & \cr
  &1      & 0 & \cdots & 1   & \cdots & 0 \\
}.
$$
Then we define the images of all $x_i$ as the following $(2n+1) \times (2n+1)$ matrices,
$$
\rho(x_i) =  
\left(
\begin{array}{ccc|cccc}
1  &        &    &  &    &        &   \\
   & \ddots &    &  &    & M_i    &   \\
   &        &  1 &  &    &        &   \\ 
   &        &    &  &    &        &   \\ \hline
   &        &    &  &  1 &        &   \\
   &        &    &  &    & \ddots &   \\
   &        &    &  &    &        & 1 \\
\end{array}
\right).
$$

Now we show that the map $\rho$ extends to an embedding of $F_n$ into $\ut{2n+1}$. Denote by $U$ the image of $F_n$. Images of all generators $x_i$ are denoted by $m_i = \rho(x_i)$. It is easy to see that for any distinct $i$ and $j$ we have $[m_i, m_j] \neq E, \ [m_i,m_j] = [m_j,m_i]^{-1}$, and $[[m_i,m_j], m_k] = E$ for any $i,j,k = \oneTo{n}$, where $E$ is the $(2n+1) \times (2n+1)$ identity matrix. Thus an image under the map $\rho$ of any word in the alphabet $X \cup X^{-1}$ can be reduced to an expression $\expForm{m}{\alpha}{n} \prod y^{\beta_{ij}}_{ij}$ in the group $U$, where $\alpha_i, \beta_{ij} \in \Z, \ i < j,  y_{ij} = [m_i, m_j]$, so the group $U$ is a two step nilpotent group with generators $\listSym{m}{n}$. To claim that the map $\rho: F_n \rightarrow UT_{2n+1}(\Z)$ is embedding, it remains to prove that the map $\rho$ has a trivial kernel. In other words, it suffices to show that $\expForm{m}{\alpha}{n} \prod y^{\beta_{ij}}_{ij} = E$ iff $\alpha_i = 0$ and $\beta_{ij} = 0, \ i,j = \oneTo{n}, i < j$. 

Let $\expForm{m}{\alpha}{n} \prod y^{\beta_{ij}}_{ij} = E$, then $\expForm{m}{\alpha}{n} = \prod y^{-\beta_{ij}}_{ij}$. Since every $y_{ij}$ commutes with $m_i, \ i = \oneTo{n},$ we  get the following,
$$
\begin{array}{l}
    [\expForm{m}{\alpha}{n}, m_i] = [\prod y^{\beta_{ij}}_{ij}, m_i],
\end{array} 
$$
$$
\begin{array}{l}
    \prod y^{\alpha_j}_{ji} = E. 
\end{array} 
$$

Recall that $U' = [U,U]$ is an abelian subgroup of $\ut{2n+1}$, so $U'$ is torsion free and the latter equality holds iff $\alpha_i = 0, \ i = \oneTo{n}$. Similarly, $\prod y^{\beta_{ij}}_{ij} = E$ iff $\beta_{ij} = 0, i,j = \oneTo{n}, \ i < j$. Therefore, $F_{\numGenerators}$ is embeddable into $\ut{53}$, so $\ut{r}, r \geq 53$, has undecidable $\KP$. Since $\ut{53}$ is a subgroup of $GL_n(\Z)$, $SL_n(\Z), \ n \geq 53$, then $GL_{n}(\Z)$, $SL_{n}(\Z)$ have undecidable $\KP$ for $n \geq 53$.

\begin{lemma}\label{lemma:KPUndecidableInQuotientGroup}
Let $G$ be a finitely generated polycyclic group and $H$ be a normal subgroup of $G$ such that the quotient group $G/H$ has undecidable $\KP$. Then group $G$ has undecidable $\KP$.
\end{lemma}
\proof Assume that the group $G$ has decidable $\KP$, that is there is an algorithm that solves $\KP$ problem in $G$. Let $A$ denote the quotient group $G/H$. Suppose we have an input for $\KP$ in the group $A$: $a_1H, a_2H, \ldots, a_kH, aH$, where $a_i, a \in G$. To solve $\KP$ we are required to find numbers $\listSym{\epsilon}{n} \in \Z$ such that
\begin{equation}\label{eq:KPinGroupA}
(a_1H)^{\epsilon_1}(a_2H)^{\epsilon_2}\ldots (a_kH)^{\epsilon_k} = a H.
\end{equation}
This equation is equivalent to the following:
$$a_1^{\epsilon_1}H a_2^{\epsilon_2}H \ldots a_k^{\epsilon_k}H = a H,$$
$$a_1^{\epsilon_1} a_2^{\epsilon_2} \ldots a_k^{\epsilon_k} H = a H,$$
$$\exists h \in H \ a_1^{\epsilon_1} a_2^{\epsilon_2} \ldots a_k^{\epsilon_k} h = a.$$

If $H$ is a finitely generated polycyclic group then there exists $b_1, \ldots, b_m \in H$ such that for any $h \in H$ there are integers $k_1, \ldots, k_m$ that $h = b_1^{k_1} \ldots b_m^{k_m}$. Hence if we solve $\KP$ problem $a_1^{\epsilon_1} a_2^{\epsilon_2} \ldots a_k^{\epsilon_k} b_1^{k_1} \ldots b_m^{k_m} = a$ in the group $G$, we get solution of $\KP$ (\ref{eq:KPinGroupA}) in the group $A$. This contradicts the assumption that the group $A$ has undecidable $\KP$. $\square$

\begin{corollary}
Let $G$ be a polycyclic group and $Fit(G)$ have rank of derived subgroup greater or equals than $322$. Then $\KP$ is undecidable in $G$.
\end{corollary}
\proof
Since $G$ is a polycylcic group then $F = Fit(G)$ is a nilpotent group. Thus $F' = F / [[F,F],F]$ is a nilpotent class two group with rank of derived subgroup greater or equal to $322$. By \textit{Theorem~\ref{theorem:main}} the $\KP$ is undecidable for $F'$ and by \textit{Lemma~\ref{lemma:KPUndecidableInQuotientGroup}} the $\KP$ is undecidable for the group $G$.

\begin{corollary}
Let $G$ be a nilpotent group of class $c \geq 3$ with lower central series
$$G = G_1 \trianglerighteq G_2 \trianglerighteq \ldots \trianglerighteq G_c \trianglerighteq G_{c+1} = \{1\},$$
where $G_{k+1} = [G_k, G]$, $k = \oneTo{c}$. Let $N$ be the quotient group $G/G_3$. If $rank([N,N]) > \numCommutators$ then the group $G$ has undecidable $\KP$.
\end{corollary}

\proof The group $N$ has undecidable $\KP$ by Corollary~\ref{corollary:freenilp}. Hence, the group $G$ has undecidable $\KP$ problem by Lemma~\ref{lemma:KPUndecidableInQuotientGroup}. $\square$

\end{document}